\documentclass[11pt,a4paper]{amsart}
\usepackage[T1]{fontenc}
\usepackage[utf8]{inputenc}
\usepackage{lmodern,mathtools,amssymb,mathrsfs,microtype}
\usepackage[a4paper,margin=25mm]{geometry}
\usepackage[hidelinks]{hyperref}
\theoremstyle{plain}
\newtheorem{maintheorem}{Theorem}

\newtheorem{theorem}{Theorem}[section]
\newtheorem{lemma}[theorem]{Lemma}
\newtheorem{proposition}[theorem]{Proposition}
\newtheorem{corollary}[theorem]{Corollary}
\theoremstyle{definition}

\theoremstyle{remark}
\newtheorem{remark}[theorem]{Remark}
\numberwithin{equation}{section}
\newcommand{\R}{\mathbb R}
\newcommand{\N}{\mathbb N}
\newcommand{\Nzero}{\mathbb N_0}
\newcommand{\B}{\mathbb B}
\newcommand{\A}{\mathbb A}
\newcommand{\ind}[1]{\mathbf 1_{#1}}
\newcommand{\norm}[1]{\lVert#1\rVert}
\DeclareMathOperator{\MAlg}{MAlg}
\DeclareMathOperator{\UC}{UC}
\DeclareMathOperator{\Lip}{Lip}
\DeclareMathOperator{\TV}{TV}
\DeclareMathOperator{\dist}{dist}
\DeclareMathOperator{\diam}{diam}
\DeclareMathOperator{\conv}{conv}
\DeclareMathOperator{\clip}{clip}
\title[Injective groups and functional balance]{Uniform retractions onto a probability measure algebra and an application to topological groups}
\author{Dekui Peng}
\address{Institute of Mathematics, Nanjing Normal University,
  Nanjing 210024, China}
\email{pengdk10@lzu.edu.cn}

\author{Li-Hong Xie}
\address{School of Mathematics and Computational Science, Wuyi University,  
Jiangmen, Guangdong 529000, P.R. China}
\email{yunli198282@126.com}

\author{Jiang Yang}
\address{School of Mathematical Sciences, Guangxi Minzu University, Nanning 530006, China}
\email{yangjiangdy@126.com}
\thanks{Jiang Yang is the corresponding author.}

\date{}
\begin{document}
\begin{abstract}

We prove that the indicator functions form a uniformly continuous
retract of the set of $[0,1]$-valued functions in $L^2([0,1])$.
Consequently, the Lebesgue probability measure algebra, equipped
with the symmetric-difference metric, is an injective uniform space.
The construction uses dyadic conditional expectations and a common
exponential coupling to obtain Lipschitz maps into the indicator
functions that approximate the identity uniformly on that set.
As an application, we construct a Polish topological group whose
left uniform space is injective but which admits a bounded left
uniformly continuous real-valued function that is not right
uniformly continuous. This answers the injective part of a question
of Bouziad and Troallic in the negative. The group is a semidirect
product of two copies of the measure algebra and is also bounded,
locally finite, and of exponent at most four.

\end{abstract}
\maketitle

\section{Introduction}\label{sec:introduction}

All topological groups and uniform spaces in this paper are Hausdorff.
Let $G$ be a topological group with identity $e$, and let $W$ range over
the neighbourhoods of $e$. The sets
\[
 L_W=\{(x,y)\in G^2:x^{-1}y\in W\},\qquad
 R_W=\{(x,y)\in G^2:xy^{-1}\in W\}
\]
form bases for the \emph{left uniformity} $\mathcal U_L$ and the
\emph{right uniformity} $\mathcal U_R$, respectively. A group is
\emph{balanced}, or a \emph{SIN group}, if $\mathcal U_L=\mathcal U_R$.
Equivalently, every identity neighbourhood contains an identity
neighbourhood invariant under conjugation. A group is \emph{functionally
balanced}, or an \emph{FSIN group}, if every bounded real-valued function
which is uniformly continuous for $\mathcal U_L$ is uniformly continuous
for $\mathcal U_R$.

Every SIN group is FSIN. The converse is related to questions of Comfort
and Ross~\cite{CR66} and Itzkowitz~\cite{Itzkowitz74}. It holds for locally
compact groups; see, in particular, Milnes~\cite{Milnes90}.
Protasov~\cite{Protasov91} proved that every almost metrizable FSIN group
is SIN. Here \emph{almost metrizable} means that the group contains a
nonempty compact subset of countable character: that subset has a
countable neighbourhood base in the group. This includes metrizable and
locally compact groups. Megrelishvili, Nickolas, and Pestov~\cite{MNP97}
obtained the corresponding result for locally connected groups.
Their theorem is stated for arbitrary real-valued functions, but its
proof uses only bounded functions; see also~\cite[Section~2]{BT07}.
Further results appear in~\cite{BT06}, and the survey~\cite{BT07}
collects questions concerning the two uniform structures.

We consider the extension question from that survey. A uniform space $Y$
is \emph{injective} if, whenever $A$ is a uniform subspace of a uniform
space $X$, every uniformly continuous map $f:A\to Y$ has a uniformly
continuous extension $\widetilde f:X\to Y$. The subspace $A$ need not be
closed. A topological group is called \emph{injective} if its left
uniform space is injective. It is called \emph{bounded} if every left
uniformly continuous real-valued function on it is bounded.

Bouziad and Troallic asked whether every injective, respectively bounded,
topological group is functionally balanced~\cite[Question~6]{BT07}.
They also observed that an injective FSIN group is SIN. Bouziad later
answered the bounded part in the negative~\cite[pp.~425--426]{Bouziad19}.
The argument uses the hereditary property of functional balance and the
Hartman--Mycielski embedding of an arbitrary topological group into a
bounded group~\cite{HM58}.

Our first main result answers the injective part. Recall that a space is
\emph{Polish} if its topology is separable and completely metrizable.
A group is \emph{locally finite} if every finitely generated subgroup is
finite. To say that its exponent divides four means that $g^4=e$ for
every element $g$.

\begin{maintheorem}\label{thm:main}
There exists an injective Polish topological group $G$ which is not
functionally balanced. The group admits a compatible complete
left-invariant metric. It is bounded, locally finite, and of exponent
at most four.
\end{maintheorem}

The main ingredient is an extension property of a measure algebra.
Let $I=[0,1]$ and let $\lambda$ be Lebesgue probability measure on $I$.
For measurable sets $A,B\subseteq I$, write
$A\triangle B=(A\setminus B)\cup(B\setminus A)$ for their symmetric
difference. The \emph{measure algebra} $\B=\MAlg(I,\lambda)$ consists of
measurable sets identified modulo null sets, with metric
$d_\lambda(A,B)=\lambda(A\triangle B)$.
For a measurable set $A$, its indicator $\ind A$ equals $1$ on $A$
and $0$ elsewhere. We identify $A\in\B$ with $\ind A$ when using
function-space notation.

Let $L^2(I;\R)$ denote the real square-integrable functions modulo
almost-everywhere equality, with norm
$\norm h_2=(\int_I|h|^2\,d\lambda)^{1/2}$. Put
\[
 C_I=\{h\in L^2(I;\R):0\le h\le1\text{ almost everywhere}\},\qquad
 K_I=\{\ind A:A\subseteq I\text{ measurable}\}.
\]
A \emph{uniform retraction} onto a uniform subspace is a uniformly
continuous map which fixes every point of that subspace.

\begin{maintheorem}\label{thm:retraction}
There is a uniformly continuous retraction from $C_I$, equipped with
its $L^2$ metric, onto $K_I$, equipped with the symmetric-difference
metric. Consequently, $(\B,d_\lambda)$ is an injective uniform space.
\end{maintheorem}

We prove Theorem~\ref{thm:retraction} first. In a product representation
of $I$, we use dyadic conditional expectations in one coordinate and
select the approximation level through the other coordinate. The common
exponential coupling of Angel and Spinka~\cite{AS21} controls the
disagreement between the selected levels for nearby inputs. A fibrewise
variance estimate gives uniform approximation on all indicators. We then
join successive approximations in shrinking neighbourhoods of the
indicator set.

Section~\ref{sec:preliminaries} fixes notation and derives the required
Hilbert-space extension property from Isbell's results~\cite{Isbell61}.
Sections~\ref{sec:approximation} and~\ref{sec:retraction} construct the
retraction and prove Theorem~\ref{thm:retraction}.
Section~\ref{sec:group} constructs the semidirect product and proves
Theorem~\ref{thm:main}.

\section{Notation and extension properties}\label{sec:preliminaries}

\subsection{Uniform spaces, metrics, and function spaces}\label{subsec:notation}

We write $\N=\{1,2,\ldots\}$, $\Nzero=\{0,1,2,\ldots\}$, and
$\R$ for the real numbers. For uniform spaces $X,Y$, the set of
uniformly continuous maps from $X$ to $Y$ is denoted by $\UC(X,Y)$.
If $X,Y$ have metrics $d_X,d_Y$, uniform continuity means that for
every $\varepsilon>0$ there is $\eta>0$ such that
$d_X(x,x')<\eta$ implies $d_Y(F(x),F(x'))<\varepsilon$ for all
$x,x'\in X$. We put
\[
 \Lip(F)=\sup_{x\ne x'}
 \frac{d_Y(F(x),F(x'))}{d_X(x,x')},
\]
allowing the value $+\infty$ and taking the supremum to be $0$ when
there are no distinct points. The map is \emph{Lipschitz} if this
number is finite. The metrics will be specified whenever this notation
is used. A \emph{uniform isomorphism} is a bijection for which both the
map and its inverse are uniformly continuous.

For a nonempty subset $M$ of a metric space $(X,d_X)$, let
\[
 \dist_{d_X}(x,M)=\inf_{y\in M}d_X(x,y),\qquad
 \diam_{d_X}(M)=\sup_{y,z\in M}d_X(y,z).
\]
We omit the metric subscript when it is fixed by the context.
For $t\in\R$, let
$\clip_{[0,1]}(t)=\min\{1,\max\{0,t\}\}$.
The symbol $\lfloor t\rfloor$ denotes the greatest integer at most $t$,
and $\{t\}=t-\lfloor t\rfloor$ is the fractional part of $t$.
Throughout, $\log$ is the natural logarithm and $\exp$ is its inverse.

All measures are completed, and measurable functions are identified
when equal almost everywhere. For a probability space $(X,\Sigma,\mu)$,
the symbol $\Sigma$ denotes its sigma-algebra of measurable sets.
For $1\le p<\infty$, we use the real space
\[
 L^p(X,\mu;\R)=\{h:X\to\R\text{ measurable}:\int_X|h|^p\,d\mu<\infty\}
 \quad\text{modulo almost-everywhere equality},
\]
with norm $\norm h_p=(\int_X|h|^p\,d\mu)^{1/p}$.
On $L^2$, the real inner product is
$\langle h,g\rangle=\int_Xhg\,d\mu$.
We suppress $\mu$ and $\R$ when the underlying space is fixed.

In a general group we use multiplicative notation and the identity $e$.
For a Boolean group, meaning an abelian group in which every element
squares to the identity, we use additive notation with identity $0$.
On a measure algebra the addition is symmetric difference and $0$ is
the empty-set class. Inversion is a uniform isomorphism from
$(G,\mathcal U_L)$ to $(G,\mathcal U_R)$. Thus injectivity and boundedness
can equivalently be defined with the right uniformity. A metric $d$ is
\emph{left invariant} if $d(zx,zy)=d(x,y)$ for all $x,y,z\in G$.
A metric is \emph{compatible} if it induces the given topology.
A compatible left-invariant metric induces $\mathcal U_L$, since
$d(x,y)=d(e,x^{-1}y)$. The maps in extension problems are maps of
uniform spaces; they need not be group homomorphisms. We refer
to~\cite{Isbell64} for general background.

\subsection{Measure algebras and the product model}\label{subsec:product}

For a probability space $(X,\Sigma,\mu)$, let
$\MAlg(X,\mu)$ be $\Sigma$ modulo the relation
$A\sim B$ if $\mu(A\triangle B)=0$. Union, intersection, and complement
descend to these equivalence classes and give the Boolean algebra
operations. Its measure metric is $d_\mu(A,B)=\mu(A\triangle B)$.
Under the indicator representation,
\begin{equation}\label{eq:indicator-metrics}
 d_\mu(A,B)=\norm{\ind A-\ind B}_1
           =\norm{\ind A-\ind B}_2^2.
\end{equation}
In particular, the $L^1$ and $L^2$ metrics on indicators induce the same
uniformity.

\begin{lemma}\label{lem:measure-algebra}
The metric $d_\mu$ on $\MAlg(X,\mu)$ is complete. The space
$(\B,d_\lambda)$ is also separable.
\end{lemma}
\begin{proof}
Let $(A_j)$ be a Cauchy sequence. Passing to a subsequence, we may
suppose that $\sum_j\mu(A_{j+1}\triangle A_j)<\infty$.
Outside a null set, membership in $A_j$ is eventually constant.
Its eventual value defines a measurable indicator, and dominated
convergence gives convergence in $d_\mu$. The original Cauchy sequence
has the same limit. For $\B$, finite unions of dyadic intervals form a
countable dense set.
\end{proof}

Let $U$ and $V$ be two copies of $I$. Their Lebesgue measures are denoted
by $\lambda_U$ and $\lambda_V$, and $\mu=\lambda_U\otimes\lambda_V$
denotes the completed product measure. We write $(u,v)$ for a point
of $U\times V$. Define
\[
 \begin{aligned}
 C&=\{h\in L^2(U\times V,\mu;\R):0\le h\le1\text{ almost everywhere}\},\\
 K&=\{\ind E:E\subseteq U\times V\text{ measurable}\}.
 \end{aligned}
\]
On $C$ we use the distance $\norm{h-g}_2$. On $K$ we use
$d_B(a,b)=\norm{a-b}_1$, where the subscript $B$ labels the Boolean
measure metric and is not a varying set. Unless another space is
specified, all subsequent $L^p$ norms refer to $(U\times V,\mu)$.
An integral over one coordinate uses $\lambda_U$ or $\lambda_V$.
For a map $F:C\to K$, its \emph{error on the indicators} is
\begin{equation}\label{eq:indicator-error}
 r(F)=\sup_{k\in K}d_B(F(k),k).
\end{equation}

The interval and product models are equivalent modulo null sets.
Represent a point of $I$ by its binary digits and separate the odd
and even positions. After discarding the null sets of exceptional
expansions, this gives a measure-preserving isomorphism between $I$
and $U\times V$. Composition with this isomorphism, also called
\emph{pullback}, preserves all $L^p$ norms and maps indicators onto
indicators. It therefore suffices to construct a uniform retraction
from $C$ onto $K$.

\subsection{The Hilbert-space extension property}\label{subsec:hilbert}

\begin{lemma}\label{lem:retract}
A uniform retract of an injective uniform space is injective.
Injectivity is preserved by uniform isomorphisms.
\end{lemma}
\begin{proof}
Let $Y\subseteq Z$, let $Z$ be injective, and let $r:Z\to Y$ be a
uniform retraction. Given a uniformly continuous map $f:A\to Y$ from
a uniform subspace $A\subseteq X$, extend its composite with the
inclusion $Y\hookrightarrow Z$ to a uniformly continuous map
$\widetilde f:X\to Z$. Then $r\circ\widetilde f$ extends $f$.
The second assertion follows by composition with a uniform isomorphism
and its inverse.
\end{proof}

\begin{lemma}\label{lem:hilbert-ball}
The closed unit ball of a real Hilbert space, with its norm uniformity,
is injective.
\end{lemma}
\begin{proof}
Let $H$ be the Hilbert space, with norm $\norm{\cdot}$ and inner product
$\langle\cdot,\cdot\rangle$, and put
$B_H=\{x\in H:\norm x\le1\}$.
Let $\mathcal C(H)$ be the family of nonempty bounded closed convex
subsets of $H$. Its Hausdorff metric is
\[
 d_H(A,D)=\max\{\sup_{a\in A}\dist(a,D),
                         \sup_{d\in D}\dist(d,A)\}.
\]
Here all point-to-set distances use the Hilbert norm.
We use two results of Isbell: $\mathcal C(H)$ is a uniform retract
of every metric space in which it is isometrically
embedded~\cite[Lemma~3.2]{Isbell61}, and the supremum-metric cube
$[0,2]^J$ is injective for every set $J$~\cite[Corollary~1.3]{Isbell61}.
The latter cube consists of functions $J\to[0,2]$, with distance
$\sup_{j\in J}|x_j-y_j|$. Isbell states the result for the unit ball
$[-1,1]^J$; coordinatewise translation gives $[0,2]^J$.
His definition of extension space allows arbitrary uniform subspaces.

\smallskip
\noindent\emph{The bounded hyperspace.}
Let
\[
 \mathscr K_H=\{A\in\mathcal C(H):A\subseteq B_H\}.
\]
Write $P_{B_H}:H\to B_H$ for metric projection: $P_{B_H}(x)$ is the
unique point of $B_H$ nearest to $x$. This map is $1$-Lipschitz.
For a subset $M\subseteq H$, let $\conv M$ denote its convex hull
and $\overline{\conv}M$ its norm-closed convex hull; the overline denotes
closure in the Hilbert norm. Taking the closed
convex hull does not increase Hausdorff distance. Indeed, finite convex
combinations with the same coefficients preserve a bound on the
distances between paired points, and that bound passes to closures.
Consequently,
\[
 T(A)=\overline{\conv}\bigl(P_{B_H}(A)\bigr)
\]
is a $1$-Lipschitz retraction from $\mathcal C(H)$ onto $\mathscr K_H$.
The space $\mathscr K_H$ is complete. To see this, use completeness of
the Hausdorff metric on nonempty closed bounded subsets of a complete
metric space. That fact follows by selecting successive points along
a subsequence with summable Hausdorff distances. Convexity and
containment in $B_H$ pass to the limit.

\smallskip
\noindent\emph{Injectivity of the bounded hyperspace.}
Suppose a metric space $X$ contains $\mathscr K_H$ isometrically.
Glue $X$ and $\mathcal C(H)$ along $\mathscr K_H$, retaining their
original metrics and defining cross distances by
\[
 d(x,A)=\inf_{D\in\mathscr K_H}
       \bigl(d_X(x,D)+d_H(D,A)\bigr),
 \qquad x\in X,\ A\in\mathcal C(H).
\]
The triangle inequality follows from agreement of the two metrics on
$\mathscr K_H$. Completeness makes the common subspace closed in both
pieces, so different points in the glued space have positive distance.
Apply Isbell's retraction onto $\mathcal C(H)$, compose with $T$, and
restrict to $X$. We obtain a uniform retraction $X\to\mathscr K_H$.

Since $\diam\mathscr K_H\le2$, the distance-coordinate map
\[
 A\longmapsto \bigl(d_H(A,D)\bigr)_{D\in\mathscr K_H}
\]
is an isometric embedding into $[0,2]^{\mathscr K_H}$.
The preceding construction gives a uniform retraction of this injective
cube onto $\mathscr K_H$. Lemma~\ref{lem:retract} shows that
$\mathscr K_H$ is injective.

\smallskip
\noindent\emph{A nearest-point selection.}
For $A\in\mathscr K_H$, let $s(A)$ be the point of $A$ nearest to the
origin. It exists and is unique by the Hilbert projection
theorem~\cite[Theorem~I.2.5]{Conway90}.
Fix $A,D\in\mathscr K_H$, and put $p=s(A)$ and $q=s(D)$.
Given $\eta>d_H(A,D)$, choose $a\in A$ and $d\in D$ with
$\norm{a-q}<\eta$ and $\norm{d-p}<\eta$.
The minimizing property gives
$\langle p,a-p\rangle\ge0$ and $\langle q,d-q\rangle\ge0$.
Adding these inequalities yields
\[
 \norm{p-q}^2
 \le \norm p\,\norm{a-q}+\norm q\,\norm{d-p}
 \le2\eta.
\]
Letting $\eta$ decrease to $d_H(A,D)$ gives
\begin{equation}\label{eq:selection}
 \norm{s(A)-s(D)}\le\sqrt{2d_H(A,D)}.
\end{equation}
The singleton map $j:B_H\to\mathscr K_H$, defined by $j(x)=\{x\}$,
is isometric and satisfies $s\circ j=\mathrm{id}_{B_H}$, where
$\mathrm{id}_{B_H}$ denotes the identity map on $B_H$.
Thus $B_H$ is a uniform retract of an injective space, and
Lemma~\ref{lem:retract} completes the proof.
\end{proof}

\begin{corollary}\label{cor:C-injective}
The uniform space $C$, equipped with its $L^2$ metric, is injective.
\end{corollary}
\begin{proof}
Since $\mu(U\times V)=1$, every $h\in C$ has $\norm h_2\le1$.
Pointwise clipping to $[0,1]$ is a $1$-Lipschitz retraction of the closed
unit ball of $L^2(U\times V,\mu;\R)$ onto $C$.
Apply Lemmas~\ref{lem:hilbert-ball} and~\ref{lem:retract}.
\end{proof}

\section{Uniform approximation of the indicator functions}
\label{sec:approximation}

Our aim is to construct Lipschitz maps from $C$ to $K$ whose error on $K$
is arbitrarily small, with error measured by $r$ in
\eqref{eq:indicator-error}.
A fixed dyadic level cannot approximate every indicator with uniform
accuracy. We therefore choose the level according to the input. We first
prepare a rounding map at each level and a common coupling of probability
distributions. These give a choice of levels for which we can control both
the Lipschitz constant and the error.

\subsection{Dyadic rounding}

For $n\in\Nzero$ and integers $0\le j<2^n$, let
\[
 J_{n,j}=[j2^{-n},(j+1)2^{-n}),
\]
with $1$ included in the last interval. Let $\mathcal D_n(V)$ be the finite
$\sigma$-algebra generated by these intervals, and let $\mathcal L(U)$ be
the Lebesgue $\sigma$-algebra on $U$. Denote by $\mathcal F_n$ the completion
of the product $\sigma$-algebra
$\mathcal L(U)\otimes\mathcal D_n(V)$ with respect to $\mu$. Thus
$(\mathcal F_n)_{n\ge0}$ is an increasing sequence of $\sigma$-algebras,
each of which retains all measurable information in the coordinate $u$.

For $h\in L^2(U\times V,\mu)$, let
$E_nh=\mathbb E_\mu[h\mid\mathcal F_n]$ be its conditional expectation
with respect to $\mathcal F_n$. This is the $\mathcal F_n$-measurable
function, unique up to a $\mu$-null set, such that
\[
 \int_H E_nh\,d\mu=\int_H h\,d\mu
 \qquad(H\in\mathcal F_n).
\]
In the present setting, it is given by the averaging formula
\begin{equation}
\label{eq:dyadic-expectation}
 E_nh(u,v)=2^n\int_{J_{n,j}}h(u,w)\,d\lambda_V(w)
 \qquad(v\in J_{n,j}).
\end{equation}
In particular, $\mathcal D_0(V)$ is the trivial $\sigma$-algebra and
$E_0h(u,v)=\int_V h(u,w)\,d\lambda_V(w)$.

The operators $E_n$ are orthogonal projections onto nested closed
subspaces of $L^2(U\times V,\mu)$. The union of their ranges contains the
finite sums of products $a(u)b(v)$, where $a\in L^2(U,\lambda_U)$ and $b$
is a dyadic step function on $V$. These sums are dense in
$L^2(U\times V,\mu)$, so $E_nh\to h$ in $L^2$. We refer to
\cite[Chapter~9]{Williams91} for these projection properties of conditional
expectation.

Using fractional parts as in Section~\ref{subsec:notation}, define
$F_n\colon C\to K$ by
\begin{equation}
\label{eq:rounding}
 F_n(h)(u,v)=\ind{\{(u,v):\{2^nv\}<E_nh(u,v)\}}(u,v).
\end{equation}
On each interval $J_{n,j}$, this replaces the mean of $h$ by the indicator
of an initial subinterval of the same relative length. In particular,
$F_n(h)$ is a measurable indicator.

\begin{lemma}
\label{lem:rounding}
For $h,g\in C$, $n\in\Nzero$ and every measurable set $A\subseteq U$,
\begin{equation}
\label{eq:fibre-lipschitz}
 \int_{A\times V}|F_n(h)-F_n(g)|\,d\mu
 =\int_{A\times V}|E_n(h-g)|\,d\mu
 \le\int_{A\times V}|h-g|\,d\mu.
\end{equation}
For $k\in K$ and $n\in\Nzero$, put
\[
 D_n^k(u)=\int_V|k(u,v)-E_nk(u,v)|^2\,d\lambda_V(v).
\]
Outside a $\lambda_U$-null set, the sequence $(D_n^k(u))_{n\ge0}$ is
decreasing, takes values in $[0,1/4]$, and satisfies
\begin{equation}
\label{eq:fibre-error}
 \int_V|F_n(k)(u,v)-k(u,v)|\,d\lambda_V(v)\le4D_n^k(u).
\end{equation}
\end{lemma}

\begin{proof}
Fix $u$ and a dyadic interval $J_{n,j}$. The two indicators in
\eqref{eq:rounding} are indicators of nested initial subintervals. Their
disagreement has measure equal to the length of $J_{n,j}$ times the
absolute difference of the two means. Summing over $j$ and integrating
over $A$ gives the equality in \eqref{eq:fibre-lipschitz}. The inequality
follows from the triangle inequality in \eqref{eq:dyadic-expectation}.

On almost every $u$-fibre, the averaging operators are nested orthogonal
projections on $L^2(V,\lambda_V)$. Hence their squared approximation
errors are decreasing. At level zero, the error is the variance of an
indicator, which is at most $1/4$.

It remains to prove \eqref{eq:fibre-error}. Let $p$ be the mean of $k$ on
$J_{n,j}$. The indicators $k$ and $F_n(k)$ have the same mean there. Their
disagreement therefore occupies at most the proportion
\[
 2\min(p,1-p)\le4p(1-p)
\]
of that interval. Its contribution to $D_n^k(u)$ is
$2^{-n}p(1-p)$. Summing over $j$ proves the assertion.
\end{proof}

The restriction to $A\times V$ in \eqref{eq:fibre-lipschitz} allows us to
compare two inputs on the fibres where they choose the same level. This
is why every $\mathcal F_n$ retains the whole coordinate $u$.

\subsection{A common exponential coupling}

Let $\ell^1(\Nzero)$ be the space of real sequences
$b=(b_n)_{n\ge0}$ for which
$\|b\|_{\ell^1}=\sum_{n\ge0}|b_n|<\infty$. The probability simplex on
$\Nzero$ is
\[
 \Delta=\left\{q=(q_n)_{n\ge0}\in\ell^1(\Nzero):
 q_n\ge0\text{ for every }n,\quad\sum_{n\ge0}q_n=1\right\}.
\]
We give $\Delta$ the $\ell^1$ metric. For $q,r\in\Delta$, their total
variation distance is
\[
 \TV(q,r)=\frac12\|q-r\|_{\ell^1}
 =\frac12\sum_{n\ge0}|q_n-r_n|.
\]

Fix independent exponential random variables $(Z_n)_{n\ge0}$ of mean one
on $(U,\lambda_U)$; equivalently,
$\lambda_U\{u:Z_n(u)>t\}=\exp(-t)$ for $t\ge0$. We use the same measurable
functions $Z_n(u)$ for all inputs. Such a family can be specified by
binary digits. Let $b_j(u)$ be the $j$th binary digit of $u$, using the
terminating expansion at dyadic points. Put
\[
 W_n(u)=\sum_{j\ge1}b_{2^n(2j-1)}(u)2^{-j},
 \qquad Z_n(u)=-\log W_n(u)
 \qquad(n\in\Nzero).
\]
The sets $\{2^n(2j-1):j\ge1\}$ partition $\N$. Hence the $W_n$ are
independent and uniformly distributed on $[0,1]$. We assign arbitrary
positive finite values to the $Z_n$ on the exceptional null sets.

For $q\in\Delta$ and $u\in U$, compare the numbers $Z_n(u)/q_n$, with
value $+\infty$ when $q_n=0$. Define $N(q,u)$ to be the unique minimizing
index when it exists, and put $N(q,u)=0$ otherwise. This defines a
measurable map $N\colon\Delta\times U\to\Nzero$: each strict-minimum
event is a countable intersection of measurable comparisons.

The following is the countable exponential coupling of Angel and Spinka
\cite[Proposition~6]{AS21}. We include the calculation used below.

\begin{lemma}
\label{lem:coupling}
For every $q\in\Delta$, the minimum defining $N(q,u)$ is uniquely
attained for $\lambda_U$-almost every $u$, and
\[
 \lambda_U\{u:N(q,u)=n\}=q_n
 \qquad(n\in\Nzero).
\]
Moreover, for $q,r\in\Delta$,
\begin{equation}
\label{eq:coupling}
 \lambda_U\{u:N(q,u)\ne N(r,u)\}
 \le\frac{2\TV(q,r)}{1+\TV(q,r)}
 \le2\TV(q,r).
\end{equation}
\end{lemma}

\begin{proof}
Fix $n$ with $q_n>0$. Conditional on $Z_n=x$, independence gives the
probability that $n$ is the strict minimum. Integrating with respect to
the density of $Z_n$, we obtain
\[
 \int_0^\infty e^{-x}\prod_{j\ne n}e^{-xq_j/q_n}\,dx=q_n.
\]
The infinite product is the limit of the probabilities for finitely many
comparisons. The strict-minimum events are disjoint, and their
probabilities sum to one. Thus the minimum is almost surely attained at
a unique index, with the asserted distribution.

Suppose that $q_nr_n>0$. Applying the same calculation to the two
vectors gives
\[
 \begin{split}
 \lambda_U\{u:N(q,u)=N(r,u)=n\}
 &=\left(\sum_{j\ge0}
      \max\left\{\frac{q_j}{q_n},\frac{r_j}{r_n}\right\}\right)^{-1}\\
 &\ge\frac{\min(q_n,r_n)}{1+\TV(q,r)}.
 \end{split}
\]
For the inequality, replace both denominators by
$\min(q_n,r_n)$ and use
$\sum_{j\ge0}\max(q_j,r_j)=1+\TV(q,r)$. Since
$\sum_{n\ge0}\min(q_n,r_n)=1-\TV(q,r)$, summing over $n$ yields
\eqref{eq:coupling}.
\end{proof}

The exceptional null set in Lemma~\ref{lem:coupling} may depend on $q$.
The prescribed value zero defines $N$ on the whole domain, and every
estimate concerns a fixed vector or a fixed pair of vectors.

\subsection{Distributions depending on the input}

For $n\in\N$, let $\Delta_n=E_n-E_{n-1}$ be the difference of the two
projections. Thus $\Delta$ denotes the probability simplex, whereas
$\Delta_n$ denotes an operator. For $h\in C$, define
\begin{equation}
\label{eq:energies}
 \begin{gathered}
 a_n(h)=\norm{\Delta_nh}_2^2\quad(n\in\N),
 \qquad S(h)=\sum_{n\ge1}a_n(h),\\
 A_0(h)=0,\qquad A_n(h)=\sum_{j=1}^n a_j(h)\quad(n\in\N).
 \end{gathered}
\end{equation}
Here $A_n(h)$ is the cumulative sum of the first $n$ terms. Orthogonality
and the convergence $E_nh\to h$ in $L^2$ give
\[
 S(h)=\norm{h-E_0h}_2^2,
 \qquad S(h)-A_n(h)=\norm{h-E_nh}_2^2.
\]
Also $S(h)\le1/4$. Indeed, on almost every $u$-fibre, put
$p=\int_Vh(u,v)\,d\lambda_V(v)$. Since $0\le h\le1$,
\[
 \int_Vh(u,v)^2\,d\lambda_V(v)-p^2
 \le p-p^2\le\frac14.
\]
Integration with respect to $\lambda_U$ gives the required bound.

Fix $\delta>0$ and $m\in\N$. For each $h\in C$, define cumulative
probabilities $Q_n(h)$ and a probability vector $q(h)=(q_n(h))_{n\ge0}$ by
\begin{equation}
\label{eq:level-distribution}
 \begin{gathered}
 Q_{-1}(h)=0,\qquad
 Q_n(h)=\left(\frac{\delta+A_n(h)}{\delta+S(h)}\right)^m
 \quad(n\in\Nzero),\\
 q_n(h)=Q_n(h)-Q_{n-1}(h)\quad(n\in\Nzero).
 \end{gathered}
\end{equation}
The sequence $Q_n(h)$ increases to one, so $q(h)\in\Delta$. These
quantities depend on $\delta$ and $m$, which we suppress in the notation.

To describe this distribution, let $p(h)=(p_n(h))_{n\ge0}\in\Delta$ be
given by
\[
 p_0(h)=\frac{\delta}{\delta+S(h)},
 \qquad p_n(h)=\frac{a_n(h)}{\delta+S(h)}\quad(n\in\N).
\]
The vector $p(h)$ depends only on $h$ and $\delta$; we suppress $\delta$
in its notation. Then $q(h)$ is
the distribution of the maximum of $m$ independent draws
from $p(h)$. Taking a maximum favours levels beyond which little variance
remains. The term $\delta$ controls how the distribution changes with
$h$.

\begin{lemma}
\label{lem:distribution}
Let $\delta>0$ and $m\in\N$, and define $q$ by
\eqref{eq:level-distribution}. For all $h,g\in C$,
\begin{equation}
\label{eq:distribution-lipschitz}
 \TV(q(h),q(g))\le\frac{m}{\delta}\norm{h-g}_2.
\end{equation}
\end{lemma}

\begin{proof}
The ranges of the operators $\Delta_n$ are mutually orthogonal. By the
difference-of-squares identity and the Cauchy--Schwarz inequality,
\[
 \begin{split}
 \sum_{n\ge1}|a_n(h)-a_n(g)|
 &\le\left(\sum_{n\ge1}\norm{\Delta_n(h-g)}_2^2\right)^{1/2}
       \left(\sum_{n\ge1}
       (\norm{\Delta_nh}_2+\norm{\Delta_ng}_2)^2\right)^{1/2}\\
 &\le\norm{h-g}_2\bigl(\sqrt{S(h)}+\sqrt{S(g)}\bigr)
 \le\norm{h-g}_2.
 \end{split}
\]
Let $b,b'$ be nonnegative summable vectors with respective sums
$s,s'\ge\delta$. Their normalized vectors satisfy
\[
 \left\|\frac{b}{s}-\frac{b'}{s'}\right\|_{\ell^1}
 \le\frac{\|b-b'\|_{\ell^1}+|s-s'|}{s}
 \le\frac2\delta\|b-b'\|_{\ell^1}.
\]
Apply this to $b=(\delta,a_1(h),a_2(h),\ldots)$ and the corresponding
vector for $g$. It follows that
\[
 \TV(p(h),p(g))\le\frac1\delta\norm{h-g}_2.
\]

The total variation distance of two $m$-fold product distributions is
at most $m$ times the distance of their factors. To see this, subtract
the products by changing one factor at a time and use the fact that a
probability vector has $\ell^1$ norm one. Taking the image under the
maximum map cannot increase total variation, since it groups terms in
the difference of the two distributions. These two observations give
\eqref{eq:distribution-lipschitz}.
\end{proof}

\subsection{Uniform approximation}

Let $\delta>0$ and $m\in\N$. Using the vector $q(h)$ from
\eqref{eq:level-distribution}, define
\begin{equation}
\label{eq:adaptive-map}
 B_{\delta,m}(h)(u,v)=F_{N(q(h),u)}(h)(u,v)
 \qquad(h\in C).
\end{equation}
For a fixed input $h$, the coordinates of $q(h)$ are numbers obtained by
integration over the whole product space. The measurable map
$u\mapsto N(q(h),u)$ then partitions $U$ according to the selected level.
Hence \eqref{eq:adaptive-map} defines a measurable indicator. It also
defines a map on equivalence classes: changing $h$ on a $\mu$-null set
leaves $q(h)$ unchanged and changes each of the countably many $F_n(h)$
only on a $\mu$-null set.

\begin{theorem}
\label{thm:approximation}
Let $\delta>0$ and $m\in\N$. The map
$B_{\delta,m}\colon(C,\norm{\cdot}_2)\to(K,d_B)$ is Lipschitz, with
\begin{equation}
\label{eq:approximation-lipschitz}
 \Lip(B_{\delta,m})\le1+\frac{2m}{\delta}.
\end{equation}
For every $t>0$, its error on the indicators satisfies
\begin{equation}
\label{eq:approximation-error}
 r(B_{\delta,m})\le4t+\exp\left(-\frac{mt}{\delta+1/4}\right).
\end{equation}
\end{theorem}

\begin{proof}
We first prove the Lipschitz estimate. Fix $h,g\in C$ and, for
$n\in\Nzero$, put
\[
 H_n=\{u\in U:N(q(h),u)=N(q(g),u)=n\}.
\]
On $H_n\times V$, the two outputs use the same dyadic level. By
\eqref{eq:fibre-lipschitz} and the disjointness of the sets $H_n$,
\[
 \sum_{n\ge0}\int_{H_n\times V}|F_n(h)-F_n(g)|\,d\mu
 \le\sum_{n\ge0}\int_{H_n\times V}|h-g|\,d\mu
 \le\norm{h-g}_1.
\]
On the remaining fibres the disagreement is at most one. By
Lemma~\ref{lem:coupling}, their measure is at most
$2\TV(q(h),q(g))$. Since $\norm{h-g}_1\le\norm{h-g}_2$,
Lemma~\ref{lem:distribution} gives \eqref{eq:approximation-lipschitz}.

We now estimate the error. Fix $k\in K$ and write $N(u)=N(q(k),u)$.
By \eqref{eq:fibre-error},
\[
 d_B(B_{\delta,m}(k),k)
 \le4\int_U D_{N(u)}^k(u)\,d\lambda_U(u).
\]
For $n\in\Nzero$, put
\[
 d_n(k)=\int_U D_n^k(u)\,d\lambda_U(u)=S(k)-A_n(k).
\]
This sequence decreases to zero. Given $t>0$, choose the least
$n\in\Nzero$ such that $d_n(k)\le t$. If $n=0$, the pointwise
monotonicity in Lemma~\ref{lem:rounding} gives
\[
 \int_U D_{N(u)}^k(u)\,d\lambda_U(u)\le d_0(k)\le t.
\]

Suppose that $n\ge1$. Then $d_{n-1}(k)>t$. The distribution of $N$ in
Lemma~\ref{lem:coupling} gives
\[
 \begin{split}
 \lambda_U\{u:N(u)<n\}
 &=Q_{n-1}(k)
 =\left(1-\frac{d_{n-1}(k)}{\delta+S(k)}\right)^m\\
 &\le\exp\left(-\frac{mt}{\delta+1/4}\right).
 \end{split}
\]
On this set, $D_{N(u)}^k(u)\le1/4$. On its complement, $N(u)\ge n$, so
$D_{N(u)}^k(u)\le D_n^k(u)$. Therefore
\[
 \int_U D_{N(u)}^k(u)\,d\lambda_U(u)
 \le t+\frac14\exp\left(-\frac{mt}{\delta+1/4}\right).
\]
Combining these estimates proves \eqref{eq:approximation-error}
uniformly in $k$.
\end{proof}

Both the selected level $N(u)$ and the fibre variance $D_n^k(u)$ may
depend on $u$. The error estimate uses pointwise monotonicity of the
variance and the measure of the event $N(u)<n$; it requires no
independence between these quantities.

\begin{corollary}
\label{cor:approximation}
For each integer $m\ge2$, put $R_m=B_{1/8,m}$. Then
\begin{equation}
\label{eq:approximation-family}
 \Lip(R_m)\le1+16m,
 \qquad r(R_m)\le\frac{1+\frac32\log m}{m}.
\end{equation}
In particular, $r(R_m)\to0$ as $m\to\infty$.
\end{corollary}

\begin{proof}
Take $\delta=1/8$ and $t=3\log m/(8m)$ in
Theorem~\ref{thm:approximation}.
\end{proof}

\begin{remark}
The map $R_9=B_{1/8,9}$ already has error less than $1/2$. Indeed,
taking $t=1/16$ in \eqref{eq:approximation-error} gives
\[
 \Lip(R_9)\le145,
 \qquad r(R_9)\le\frac14+\exp(-3/2)<\frac12.
\]
The retraction constructed in the next section uses the convergence of
the whole family in \eqref{eq:approximation-family}.
\end{remark}

\section{A uniformly continuous retraction}\label{sec:retraction}

Corollary~\ref{cor:approximation} gives uniformly continuous maps
$C\to K$ with arbitrarily small error on $K$. We now use these maps
in shrinking neighbourhoods of $K$. A joining operation on indicators
will make the successive changes uniformly summable.

For $a,b\in K$ and $s\in[0,1]$, define $P_s:K\times K\to K$ by
\begin{equation}\label{eq:joining}
 P_s(a,b)(u,v)=
 \begin{cases}
  b(u,v),&0\le u\le s,\\
  a(u,v),&s<u\le1.
 \end{cases}
\end{equation}
Thus $P_s$ replaces $a$ by $b$ on the strip $[0,s]\times V$.
As equivalence classes, $P_0(a,b)=a$ and $P_1(a,b)=b$.
For $a,b,a',b'\in K$ and $s,t\in[0,1]$, we have
\begin{align}
 d_B(P_s(a,b),P_t(a',b'))
 &\le d_B(a,a')+d_B(b,b')+|s-t|,\label{eq:joining-continuity}\\
 d_B(P_s(a,b),a)&\le d_B(a,b).\label{eq:joining-change}
\end{align}
To obtain the first inequality, compare the two outputs outside the
strip between $s$ and $t$ and bound that strip's measure by $|s-t|$.
The second holds because a change is confined to the disagreement set
of $a$ and $b$.

\begin{lemma}\label{lem:gluing}
Suppose that for every $\varepsilon>0$ there exists $T\in\UC(C,K)$
with $r(T)\le\varepsilon$. Then there is a uniformly continuous
retraction $R:C\to K$.
\end{lemma}
\begin{proof}
We construct a sequence of uniformly continuous maps whose increments
are summable in the uniform metric.

\smallskip
\noindent\emph{Choice of maps and neighbourhoods.}
For $n\in\Nzero$, put $\varepsilon_n=2^{-n-4}$ and choose
$T_n\in\UC(C,K)$ with $r(T_n)\le\varepsilon_n$.
Having fixed all the maps $T_n$, choose positive numbers $\rho_n$
recursively so that $\rho_{n+1}<\rho_n/2$ and
\begin{equation}\label{eq:radii}
 \norm{h-g}_2\le2\rho_n
 \quad\Longrightarrow\quad
 d_B(T_i(h),T_i(g))\le\varepsilon_n
 \quad(i=n,n+1).
\end{equation}
At the $n$th stage this requires only uniform continuity of $T_n$
and $T_{n+1}$, together with a further decrease of $\rho_n$ if $n\ge1$.

For $h\in C$, let
\[
 d(h)=\dist_{L^2}(h,K)=\inf_{k\in K}\norm{h-k}_2.
\]
This is a $1$-Lipschitz real-valued function on $C$.
Define the cutoff $\chi_n:C\to[0,1]$ by
\begin{equation}\label{eq:cutoff}
 \chi_n(h)=\clip_{[0,1]}
       \left(\frac{\rho_n-d(h)}{\rho_n-\rho_{n+1}}\right).
\end{equation}
It is Lipschitz, equals $0$ when $d(h)\ge\rho_n$, and equals $1$
when $d(h)\le\rho_{n+1}$.

\smallskip
\noindent\emph{Successive joining.}
Define maps $S_n:C\to K$ by
\begin{equation}\label{eq:successive-joining}
 S_0=T_0,\qquad
 S_{n+1}(h)=P_{\chi_n(h)}(S_n(h),T_{n+1}(h)).
\end{equation}
Every $S_n$ is uniformly continuous by~\eqref{eq:joining-continuity}.
These $K$-valued maps are distinct from the scalar energy $S(h)$
defined in Section~\ref{sec:approximation}.
The endpoint property of $P_s$ gives
\begin{equation}\label{eq:near-indicators}
 S_n(h)=T_n(h)\quad\text{whenever }d(h)\le\rho_n.
\end{equation}
Indeed, the assertion holds for $n=0$. At the next stage,
$d(h)\le\rho_{n+1}$ implies $\chi_n(h)=1$, hence
$S_{n+1}(h)=T_{n+1}(h)$.

There is no change from $S_n(h)$ to $S_{n+1}(h)$ when
$d(h)\ge\rho_n$. If $d(h)<\rho_n$, then~\eqref{eq:near-indicators}
applies and we can choose $k\in K$ with $\norm{h-k}_2<2\rho_n$.
By~\eqref{eq:radii},
\begin{align*}
 d_B(T_n(h),T_{n+1}(h))
 &\le d_B(T_n(h),T_n(k))+d_B(T_n(k),k)\\
 &\quad+d_B(k,T_{n+1}(k))+d_B(T_{n+1}(k),T_{n+1}(h))\\
 &\le2\varepsilon_n+r(T_n)+r(T_{n+1})
 \le4\varepsilon_n.
\end{align*}
Combining this with~\eqref{eq:joining-change}, we obtain the global estimate
\begin{equation}\label{eq:summable-changes}
 \sup_{h\in C}d_B(S_{n+1}(h),S_n(h))\le4\varepsilon_n.
\end{equation}

\smallskip
\noindent\emph{The limit map.}
Since $\sum_n\varepsilon_n<\infty$ and $K$ is complete by
Lemma~\ref{lem:measure-algebra}, the sequence $(S_n)$ converges uniformly
to a map $R:C\to K$.
To check uniform continuity, choose a finite stage with sufficiently
small uniform error and then use the uniform continuity of that stage.
For each $k\in K$, equation~\eqref{eq:near-indicators} gives
$S_n(k)=T_n(k)$. Thus $d_B(S_n(k),k)\le\varepsilon_n\to0$,
and $R(k)=k$.
\end{proof}

\begin{proof}[Proof of Theorem~\ref{thm:retraction}]
Corollary~\ref{cor:approximation} supplies the maps required by
Lemma~\ref{lem:gluing}, so there is a uniform retraction $R:C\to K$.
The measure-preserving isomorphism in Section~\ref{subsec:product}
induces an isometry $J:C_I\to C$ which maps $K_I$ onto $K$.
Hence $J^{-1}\circ R\circ J:C_I\to K_I$ is the required retraction
on the interval model.

By Corollary~\ref{cor:C-injective}, $C$ is injective.
Equation~\eqref{eq:indicator-metrics} identifies the uniformity of
$(K,d_B)$ with its subspace uniformity from $C$.
Lemma~\ref{lem:retract} therefore makes $K$, and hence
$(\B,d_\lambda)$, injective.
\end{proof}

\section{The topological group}\label{sec:group}

We now prove Theorem~\ref{thm:main}. We first construct a continuous action
of one measure algebra on another. The resulting semidirect product has a left
uniform space isomorphic to $(\B,d_\lambda)$. We then use parity events to
construct a bounded left uniformly continuous function which is not right
uniformly continuous.

\subsection{A continuous action on a measure algebra}\label{subsec:action}

Equip $I=[0,1]$ with its Lebesgue $\sigma$-algebra and $\{0,1\}$ with its
discrete $\sigma$-algebra. For $b\in\{0,1\}$, let $\delta_b$ denote the unit
point mass at $b$. Put
\[
  \Omega=(I\times\{0,1\})^{\N},
  \qquad
  \nu=\left(\lambda\otimes\frac12(\delta_0+\delta_1)\right)^{\N}.
\]
We take the completion of this product measure, and denote its completed
$\sigma$-algebra by $\Sigma$. Let
\[
  \A=\MAlg(\Omega,\nu),
  \qquad
  d_\nu(E,F)=\nu(E\triangle F)\quad(E,F\in\A).
\]
For each $i\in\N$, let $t_i:\Omega\to I$ and
$\epsilon_i:\Omega\to\{0,1\}$ be the coordinate projections. Thus a point of
$\Omega$ has the form
\[
  \omega=\bigl((t_i(\omega),\epsilon_i(\omega))\bigr)_{i\in\N}.
\]
We give both $\A$ and $\B$ their Boolean group structures, with symmetric
difference as the group operation. In either measure algebra, $0$ denotes
the class of the empty set.

For $A\in\B$, choose a measurable representative and define
$T_A:\Omega\to\Omega$ by
\begin{equation}\label{eq:action}
  T_A(\omega)
  =\bigl((t_i(\omega),
       \epsilon_i(\omega)+\ind{A}(t_i(\omega))\pmod 2)\bigr)_{i\in\N}.
\end{equation}
On the $i$th coordinate, this map interchanges the two equally weighted bits
when $t_i(\omega)\in A$. Hence $T_A$ is measure preserving and
$T_A^2$ is the identity on $\Omega$. Its pullback induces the map
\[
  \alpha_A:\A\longrightarrow\A,
  \qquad
  \alpha_A(E)=T_A^{-1}(E),
\]
where the right-hand side denotes the class of the inverse image of a
measurable representative of $E$. Since $T_A$ preserves $\nu$, this class
does not depend on the representative of $E$. If two representatives of
$A$ differ on a null set, the associated transformations agree outside the
null event that some $t_i$ belongs to that set. Thus $\alpha_A$ also does
not depend on the representative of $A$.

Inverse images preserve the Boolean operations, and $T_A$ is a
measure-preserving involution. Therefore $\alpha_A$ is an isometric Boolean
automorphism of $\A$. For $A,B\in\B$, formula~\eqref{eq:action} gives
\[
  \alpha_A\alpha_B=\alpha_{A\triangle B}.
\]
We denote this action of $\B$ on $\A$ by $\alpha$.

\begin{lemma}\label{lem:action}
The map
\[
  \B\times\A\longrightarrow\A,
  \qquad (A,E)\longmapsto\alpha_AE,
\]
is jointly continuous for the measure metrics.
\end{lemma}

\begin{proof}
Let $s\in\N$, and suppose that $E\in\A$ has a representative depending
only on the first $s$ coordinate pairs. Membership in $E$ can differ under
$T_A$ and $T_B$ only if at least one of $t_1,\ldots,t_s$ belongs to
$A\triangle B$. Consequently,
\begin{equation}\label{eq:cylinder-estimate}
  d_\nu(\alpha_AE,\alpha_BE)\le s\,d_\lambda(A,B)
  \qquad(A,B\in\B).
\end{equation}
Events depending on finitely many coordinate pairs are dense in $\A$.
Indeed, they generate the product $\sigma$-algebra, and approximation in
measure passes through monotone limits and completion.

For an arbitrary $E\in\A$, choose an approximating finite-coordinate event
$E'$. The isometry of each $\alpha_A$ gives
\[
  d_\nu(\alpha_AE,\alpha_BE)
  \le 2d_\nu(E,E')+d_\nu(\alpha_AE',\alpha_BE').
\]
Together with~\eqref{eq:cylinder-estimate}, this proves continuity in $A$
for each fixed $E$. Finally, for $A,B\in\B$ and $E,F\in\A$, we have
\[
  d_\nu(\alpha_AE,\alpha_BF)
  \le d_\nu(E,F)+d_\nu(\alpha_AF,\alpha_BF).
\]
This proves joint continuity.
\end{proof}

Form the semidirect product $G=\A\rtimes_\alpha\B$. Its underlying set is
$\A\times\B$, and its multiplication is
\begin{equation}\label{eq:group-law}
  (E,A)(F,B)=(E\triangle\alpha_AF,A\triangle B)
  \qquad(E,F\in\A,\ A,B\in\B).
\end{equation}
The identity of $G$ is $e=(0,0)$, and
\[
  (E,A)^{-1}=(\alpha_AE,A).
\]
We equip $G$ with the product of the measure-metric topologies on $\A$ and
$\B$. By Lemma~\ref{lem:action}, multiplication and inversion are
continuous. Hence $G$ is a topological group.

\begin{proposition}\label{prop:metric}
The metric $D$ on $G$ defined by
\begin{equation}\label{eq:group-metric}
  D\bigl((E,A),(F,B)\bigr)=d_\nu(E,F)+d_\lambda(A,B)
\end{equation}
is compatible, complete and left invariant. Moreover, $(G,D/2)$ is
isometric to $(\B,d_\lambda)$. In particular, $G$ is Polish and its left
uniform space is injective.
\end{proposition}

\begin{proof}
The metric $D$ induces the product topology by definition. To prove left
invariance, fix $(H,C)\in G$. The difference of the first coordinates after
left multiplication is
\[
  (H\triangle\alpha_CE)\triangle(H\triangle\alpha_CF)
  =\alpha_C(E\triangle F).
\]
Its measure is $\nu(E\triangle F)$, since $\alpha_C$ is an isometry.
The difference of the second coordinates remains $A\triangle B$.
Thus left multiplication preserves $D$.

We next identify $(\A,d_\nu)$ with $(\B,d_\lambda)$. Expand each coordinate
$t_i$ in binary, and enumerate all these digits together with the bits
$\epsilon_i$ as one sequence. These are independent fair bits. The binary
expansion map therefore gives a measure-preserving isomorphism from
$(\Omega,\Sigma,\nu)$ to $(I,\lambda)$ modulo null sets. Let
\[
  \varphi:(\A,d_\nu)\longrightarrow(\B,d_\lambda)
\]
be the induced isometry.

For a measurable set $C\subseteq I$, write
\[
  \tfrac12 C=\{t/2:t\in C\},
  \qquad
  \tfrac12+\tfrac12 C=\{(1+t)/2:t\in C\}.
\]
These operations also define maps on measure-algebra classes. Define
\[
  \Phi:G\longrightarrow\B,
  \qquad
  \Phi(E,A)=\tfrac12\varphi(E)
             \cup\bigl(\tfrac12+\tfrac12 A\bigr).
\]
Splitting a measurable subset of $I$ into its intersections with the two
half-intervals gives the inverse of $\Phi$. The two halves meet only on a
null set, and scaling either half multiplies measure by $1/2$. Hence
\[
  d_\lambda\bigl(\Phi(E,A),\Phi(F,B)\bigr)
  =\tfrac12 D\bigl((E,A),(F,B)\bigr).
\]
This proves the asserted isometry. By Lemma~\ref{lem:measure-algebra},
$(\B,d_\lambda)$ is complete and separable. Therefore $D$ is complete
and $G$ is Polish. Since $D$ is a compatible left-invariant metric, it
induces the left uniformity of $G$. Theorem~\ref{thm:retraction} and
Lemma~\ref{lem:retract} now give injectivity of this uniform space.
\end{proof}

\subsection{A function witnessing the failure of functional balance}
\label{subsec:nonbalance}

For each $n\in\N$, define the interval $J_n\in\B$ and the parity event
$P_n\in\A$ by
\[
  J_n=[0,1/(2n)],
  \qquad
  P_n=\left\{\omega\in\Omega:
       \sum_{i=1}^n\epsilon_i(\omega)\equiv0\pmod2\right\}.
\]
Here $J_n$ denotes a shrinking interval; the dyadic intervals introduced
earlier are denoted by $J_{n,j}$ with two indices.

Let $A\in\B$ and put $p=\lambda(A)\in[0,1]$. Membership in $P_n$ changes under
$T_A$ precisely when an odd number of $t_1,\ldots,t_n$ belong to $A$.
The corresponding indicators are independent Bernoulli random variables
with parameter $p$. Since
\[
  \int_\Omega
     (-1)^{\sum_{i=1}^n\ind{A}(t_i(\omega))}\,d\nu(\omega)
  =(1-2p)^n,
\]
we obtain
\begin{equation}\label{eq:parity-distance}
  d_\nu(P_n,\alpha_AP_n)=\frac{1-(1-2p)^n}{2}.
\end{equation}
Taking $A=J_n$ gives
\begin{equation}\label{eq:parity-lower}
  d_\nu(P_n,\alpha_{J_n}P_n)
  =\frac{1-(1-1/n)^n}{2}
  \ge\frac{1-\exp(-1)}{2}>\frac14.
\end{equation}

For distinct $j,n\in\N$, we also have
\begin{equation}\label{eq:parity-separated}
  d_\nu(\alpha_{J_n}P_n,P_j)=\frac12.
\end{equation}
To verify this equality, condition on all the coordinates $t_i$.
The difference of the two parity sums contains the nonempty block of fair
bits with indices $\min(j,n)<i\le\max(j,n)$. Its sum modulo two remains
fair after all other bits have been fixed. The terms contributed by the
$t_i$ are constants under this conditioning. Integrating proves
\eqref{eq:parity-separated}.

\begin{proposition}\label{prop:witness}
The function $f:G\to[0,1]$ defined by
\begin{equation}\label{eq:witness}
  f(E,A)=\min\left\{1,\,4\inf_{j\in\N}d_\nu(E,P_j)\right\}
\end{equation}
is left uniformly continuous and is not right uniformly continuous.
\end{proposition}

\begin{proof}
Distance to a nonempty subset of a metric space is a $1$-Lipschitz
function. Thus $f$ is $4$-Lipschitz for $D$ and takes values in $[0,1]$.
By Proposition~\ref{prop:metric}, it is left uniformly continuous.

For $n\in\N$, put
\[
  g_n=(P_n,0),\qquad
  u_n=(0,J_n),\qquad
  h_n=u_ng_n=(\alpha_{J_n}P_n,J_n).
\]
Since $\lambda(J_n)=1/(2n)\to0$, we have
\[
  h_ng_n^{-1}=u_n\longrightarrow e.
\]
On the other hand, $f(g_n)=0$, while
\eqref{eq:parity-lower} and~\eqref{eq:parity-separated} imply $f(h_n)=1$.
Right uniform continuity would force
$|f(h_n)-f(g_n)|\to0$. This contradiction proves the assertion.
\end{proof}

\subsection{Boundedness and the algebraic properties}\label{subsec:algebra}

We first verify the boundedness needed in Theorem~\ref{thm:main}.

\begin{lemma}\label{lem:bounded}
Every uniformly continuous function from $(\B,d_\lambda)$ to $\R$ is bounded.
\end{lemma}

\begin{proof}
Let $v:\B\to\R$ be uniformly continuous. Choose $\eta>0$ such that
\[
  d_\lambda(A,B)<\eta\quad\Longrightarrow\quad|v(A)-v(B)|<1
  \qquad(A,B\in\B).
\]
Take an integer $N\ge1$ with $1/N<\eta$, and partition $I$ into measurable
intervals $I_1,\ldots,I_N$ of measure $1/N$. For $A\in\B$, put $C_0=0$
and
\[
  C_k=A\cap\bigcup_{i=1}^k I_i\qquad(1\le k\le N).
\]
Then $C_N=A$ and $d_\lambda(C_{k-1},C_k)\le1/N<\eta$ for every
$1\le k\le N$. Consequently,
\[
  |v(A)-v(0)|
  \le\sum_{k=1}^N|v(C_k)-v(C_{k-1})|<N.
\]
The integer $N$ is independent of $A$, so $v$ is bounded.
\end{proof}

\begin{proof}[Proof of Theorem~\ref{thm:main}]
Let $G$ be the group defined by~\eqref{eq:group-law}.
Proposition~\ref{prop:metric} shows that $G$ is Polish and injective, and
that it admits a compatible complete left-invariant metric.
By Proposition~\ref{prop:witness}, $G$ is not functionally balanced.
The uniform isomorphism in Proposition~\ref{prop:metric} and
Lemma~\ref{lem:bounded} show that every left uniformly continuous
real-valued function on $G$ is bounded.

It remains to verify the algebraic assertions. For $(E,A)\in G$,
formula~\eqref{eq:group-law} gives
\[
  (E,A)^2=(E\triangle\alpha_AE,0).
\]
The subgroup $\A\times\{0\}$ is Boolean. Hence every fourth power in $G$
is the identity, so $G$ has exponent at most four.

Finally, let $(E_1,A_1),\ldots,(E_r,A_r)$ be finitely many elements of $G$,
where $r\in\N$. Let $Q$ be the subgroup of $\B$ generated by
$A_1,\ldots,A_r$. This is a finite Boolean group. Let $H$ be the subgroup
of $\A$ generated by the finite family
\[
  \{\alpha_AE_i:A\in Q,\ 1\le i\le r\}.
\]
Then $H$ is finite. Moreover, for $A,B\in Q$ and $1\le i\le r$,
\[
  \alpha_B(\alpha_AE_i)=\alpha_{B\triangle A}E_i\in H,
\]
so $H$ is invariant under the action of $Q$. Therefore $H\rtimes_\alpha Q$
is a finite subgroup of $G$ containing all the given elements. Every
finitely generated subgroup of $G$ is thus finite. This completes the
proof of Theorem~\ref{thm:main}.
\end{proof}

\end{document}